\documentclass{article}
\usepackage{amsmath,amssymb,amsthm}
\theoremstyle{definition}
\newtheorem*{acknowledgement}{Acknowledgement}
\newtheorem{definition}{Definition}
\newtheorem{corollary}[definition]{Corollary}

\newtheorem{problem}{Problem}
\theoremstyle{plain}
\newtheorem{theorem}[definition]{Theorem}
\newtheorem{proposition}[definition]{Proposition}
\newtheorem{assertion}[definition]{Assertion}
\newtheorem{conjecture}[definition]{Conjecture}
\theoremstyle{remark}
\newtheorem{remark}[definition]{Remark}

\author{{\bf Steven Duplij} \thanks{E-mails:
{\tt Steven.A.Duplij@univer.kharkov.ua}
and {\tt duplij@member.ams.org}} \thanks{Internet:
{\tt http://gluon.physik.uni-kl.de/\~{}duplij}}\\
Theory Group,  Nuclear Physics Laboratory\\
Kharkov National University, Kharkov 61001, Ukraine
}
\title{\bf On supermatrix idempotent operator semigroups}
\begin{document}
\maketitle
\begin{abstract}
One-parameter semigroups of antitriangle idempotent supermatrices and
corresponding superoperator semigroups are introduced and
investigated. It is shown that $t$-linear idempotent superoperators
and exponential superoperators are mutually dual in some sense, and
the first gives additional to exponential solution to the
initial Cauchy problem. The corresponding functional equation and
analog of resolvent are found for them. Differential and functional
equations for idempotent (super)operators are derived for their
general $t$ power-type dependence. \end{abstract}
\newpage
\mbox{}
\vskip 5em
\section{Introduction}

Operator semigroups \cite{hil/phi} play an important role in mathematical
physics \cite{davies,goldstein,hille} viewed as a general theory of evolution
systems \cite{eng/nag,belleni,dan/koc}. Its development covers many new fields
\cite{ber/chr/res,satyan,ber/jun/mil,ahmed}, but one of vital for modern
theoretical physics directions --- supersymmetry and related mathematical
structures --- was not considered before in application to operator semigroup
theory. The main difference between previous considerations is the fact that
among building blocks (e.g. elements of corresponding matrices) there exist
noninvertible objects (divisors of zero and nilpotents) which by themselves
can form another semigroup. Therefore, we have to take that into account and
investigate it properly, which can be called a $semigroup\times semigroup$
method.

Here we study continuous supermatrix representations of idempotent operator
semigroups firstly introduced in \cite{dup12,dup15} for bands. Usually matrix
semigroups are defined over a field $\mathbb{K}$ \cite{okninski1} (on some
nonsupersymmetric generalizations of $\mathbb{K}$-representations see
\cite{pon2,okn/pon}). But after discovery of supersymmetry
\cite{vol/aku0,wes/zum1} the realistic unified particle theories began to be
considered in superspace \cite{sal/str1,gat/gri/roc/sie}. So all variables and
functions were defined not over a field $\mathbb{K}$, but over
Grassmann-Banach superalgebras over $\mathbb{K}$
\cite{dewitt,khrennikov,vla/vol}, becoming in general noninvertible, and
therefore they should be considered by semigroup theory, which was claimed in
\cite{dup6,dup10}, and some semigroups having nontrivial abstract properties
were found \cite{dup11}. Also, it was shown that supermatrices of the special
(antitriangle) shape can form various strange and sandwich semigroups not
known before \cite{dup12,dup-hab}. Here we consider one-parametric semigroups
(for general theory see \cite{davies,eng/nag,cle/hei/ang}) of antitriangle
supermatrices and corresponding superoperator semigroups. The first ones
continuously represent idempotent semigroups and second ones lead to new
superoperator semigroups with nontrivial properties.

Let $\Lambda$ be a commutative Banach $\mathbb{Z}_{2}$-graded superalgebra
\cite{berezin,kac} over a field $\mathbb{K}$ (where $\mathbb{K}=\mathbb{R},$
$\mathbb{C}$ or $\mathbb{Q}_{p}$) with a decomposition into the direct sum:
$\Lambda=\Lambda_{0}\oplus\Lambda_{1}$. The elements $a$ from $\Lambda_{0}$
and $\Lambda_{1}$ are homogeneous and have the fixed even and odd parity
defined as $\left|  a\right|  \overset{def}{=}\left\{  i\in\left\{
0,1\right\}  =\mathbb{Z}_{2}|\,a\in\Lambda_{i}\right\}  $. The even
homomorphism $\frak{r}_{body}:\Lambda\rightarrow\mathbb{B}$ is called a body
map and the odd homomorphism $\frak{r}_{soul}:\Lambda\rightarrow\mathbb{S}$ is
called a soul map \cite{rog1}, where $\mathbb{B}$ and $\mathbb{S}$ are purely
even and odd algebras over $\mathbb{K}$ and $\Lambda=\mathbb{B}\oplus
\mathbb{S}$. It can be thought that, if we have the Grassmann algebra
$\Lambda$ with generators $\xi_{i},\ldots,\xi_{n}$ $\xi_{i}\xi_{j}+\xi_{j}%
\xi_{i}=0,\,1\leq i,\,j\leq n,$ in particular $\xi_{i}^{2}=0$ ($n$ can be
infinite, and only this case is nontrivial and interesting), then any even $x$
and odd $\varkappa$ elements have the expansions (which can be infinite)
\begin{align}
x  &  =x_{body}+x_{soul}=x_{body}+x_{12}\xi_{1}\xi_{2}+x_{13}\xi_{1}\xi
_{3}+\ldots=\nonumber\\
&  x_{body}+\mathrel{\mathop{\sum}\limits_{1\leq r\leq n}}\,\mathrel
{\mathop{\sum}\limits_{1<i_{1}<\ldots<i_{2r}\leq n}}x_{i_{1}\ldots i_{2r}}%
\xi_{i_{1}}\ldots\xi_{i_{2r}}\label{x1}\\
\varkappa &  =\varkappa_{soul}=\varkappa_{1}\xi_{1}+\varkappa_{2}\xi
_{2}+\ldots+x_{123}\xi_{1}\xi_{2}\xi_{3}+\ldots=\nonumber\\
&  \mathrel{\mathop{\sum}\limits_{1\leq r\leq n}}\,\mathrel{\mathop{\sum
}\limits_{1<i_{1}<\ldots<i_{r}\leq n}}\varkappa_{i_{1}\ldots i_{2r-1}}%
\xi_{i_{1}}\ldots\xi_{i_{2r-1}} \label{x2}%
\end{align}

From (\ref{x1})-(\ref{x2}) it follows

\begin{corollary}
\label{cor}The equations $x^{2}=0$ and $x\varkappa=0$ can have nonzero
nontrivial solutions (appearing zero divisors and even nilpotents, while odd
objects are always nilpotent).
\end{corollary}

\begin{conjecture}
If zero divisors and nilpotents will be included in the following analysis as
elements of matrices, then one can find new and unusual properties of
corresponding semigroups.
\end{conjecture}

For that we should consider general properties of supermatrices \cite{berezin}
and introduce their additional reductions \cite{dup12}.

\section{Supermatrices and their even-odd classification}

Let us consider $\left(  p|q\right)  $-dimensional linear model superspace
$\Lambda^{p|q}$ over $\Lambda$ (in the sense of \cite{berezin,lei1}) as the
even sector of the direct product $\Lambda^{p|q}=\Lambda_{0}^{p}\times
\Lambda_{1}^{q}$ \cite{rog1,vla/vol}. The even morphisms $\mbox{Hom}%
_{0}\left(  \Lambda^{p|q},\Lambda^{m|n}\right)  $ between superlinear spaces
$\Lambda^{p|q}\rightarrow\Lambda^{m|n}$ are described by means of $\left(
m+n\right)  \times\left(  p+q\right)  $-supermatrices \cite{berezin,lei1} (for
some nontrivial properties see \cite{bac/fel1,urr/mor1}). In what follows we
will treat noninvertible morphisms \cite{nashed,dav/rob} on a par with
invertible ones \cite{dup12}.

First we consider $\left(  1+1\right)  \times\left(  1+1\right)
$-supermatrices describing the elements from $\mbox{Hom}_{0}\left(
\Lambda^{1|1},\Lambda^{1|1}\right)  $ in the standard $\Lambda^{1|1}$ basis
\cite{berezin}
\begin{equation}
M\equiv\left(
\begin{array}
[c]{cc}%
a & \alpha\\
\beta &  b
\end{array}
\right)  \in\mbox{{\rm Mat}}_{\Lambda}\left(  1|1\right)  \label{1}%
\end{equation}

\noindent where $a,b\in\Lambda_{0},\,\alpha,\beta\in\Lambda_{1},$ $\alpha
^{2}=\beta^{2}=0$ (in the following we use Latin letters for elements from
$\Lambda_{0}$ and Greek letters for ones from $\Lambda_{1}$, and all odd
elements are nilpotent of index 2).

The supertrace and Berezinian (superdeterminant) are defined by
\cite{berezin}
\begin{equation}
\mbox{str}M=a-b, \label{2}%
\end{equation}%
\begin{equation}
\mbox{Ber}M=\frac{a}{b}+\frac{\beta\alpha}{b^{2}}. \label{3}%
\end{equation}

First term corresponds to triangle supermatrices, second term - to
antitriangle ones. So we obviously have different \textbf{two dual types} of
supermatrices \cite{dup12}.

\begin{definition}
\label{even-odd}\textsl{Even-reduced supermatrices} are elements from
$\mbox{{\rm
Mat}}_{\Lambda}\left(  1|1\right)  $ of the form
\begin{equation}
M_{even}\equiv\left(
\begin{array}
[c]{cc}%
a & \alpha\\
0 & b
\end{array}
\right)  \in\mbox{{\rm RMat}}_{\Lambda}^{\,even}\left(  1|1\right)
\subset\mbox{{\rm Mat}}_{\Lambda}\left(  1|1\right)  . \label{4}%
\end{equation}

\noindent\textsl{Odd-reduced supermatrices} are elements from $\mbox{{\rm
Mat}}_{\Lambda}\left(  1|1\right)  $ of the form
\begin{equation}
M_{odd}\equiv\left(
\begin{array}
[c]{cc}%
0 & \alpha\\
\beta &  b
\end{array}
\right)  \in\mbox{{\rm RMat}}_{\Lambda}^{\,odd}\left(  1|1\right)
\subset\mbox{{\rm Mat}}_{\Lambda}\left(  1|1\right)  . \label{5}%
\end{equation}
\end{definition}

\begin{conjecture}
\label{b=0}The odd-reduced supermatrices have a nilpotent (but nonzero)
Berezinian
\begin{equation}
\mbox{{\rm Ber}}M_{odd}=\dfrac{\beta\alpha}{b^{2}}\neq0,\;\;\left(  \mbox{{\rm
Ber}}M_{odd}\right)  ^{2}=0. \label{b0}%
\end{equation}
\end{conjecture}

\begin{remark}
Indeed this property (\ref{b0}) prevented in the past the use of this type
(odd-reduced) of supermatrices in physics. All previous applications
(excluding \cite{dup12,dup15,dup13}) were connected with triangle
(even-reduced, similar to Borel) ones and first term in Berezinian $\mbox{{\rm
Ber}}M=\frac{a}{b}$ (\ref{3}).
\end{remark}

The even- and odd-reduced supermatrices are \textit{mutually dual} in the
sense of the Berezinian addition formula \cite{dup12}
\begin{equation}
\mbox{Ber}M=\mbox{Ber}M_{even}+\mbox{Ber}M_{odd}. \label{7a}%
\end{equation}

For sets of matrices we use corresponding bold symbols,
and the set product is standard
\[
\mathbf{M\cdot N}\overset{def}{=}\left\{  \cup MN\,|\,M,N\in\mbox{{\rm Mat}%
}_{\Lambda}\left(  1|1\right)  \right\}  .
\]
The matrices from $\mbox{Mat}\left(  1|1\right)  $ form a linear semigroup of
$\left(  1+1\right)  \times\left(  1+1\right)  $-super\-mat\-rices under the
standard supermatrix multiplication $\frak{M}\left(  1|1\right)  \overset
{def}{=}\left\{  \mathbf{M\,}|\,\cdot\right\}  $ \cite{berezin}. Obviously,
the even-reduced matrices $\mathbf{M}_{even}$ form a semigroup $\frak{M}%
_{even}\left(  1|1\right)  $ which is a subsemigroup of $\frak{M}\left(
1|1\right)  $, because of $\mathbf{M}_{even}\mathbf{\cdot M}_{even}%
\mathbf{\subseteq M}_{even}$ and the unity is in $\frak{M}_{even}\left(
1|1\right)  $. This trivial observation leads to general structure (Borel)
theory for matrices: triangle matrices form corresponding substructures
(subgroups and subsemigroups). It was believed before that in case of
supermatrices the situation does not changed, because supermatrix
multiplication is the same \cite{berezin}. But they did not take into account
\textit{zero divisors and nilpotents} appearing naturally and inevitably in supercase.

\begin{conjecture}
Standard (lower/upper) triangle supermatrices are not the only substructures
due to unusual properties of \textit{zero divisors and nilpotents }appearing
among elements (see (\ref{x1})-(\ref{x2}) and \textbf{Corollary \ref{cor}}).
\end{conjecture}

It means that in such consideration we have additional (to triangle) class of
subsemigroups. Then we can formulate the following general

\begin{problem}
For a given $n,m,p,q$ to describe and classify all possible substructures
(subgroups and subsemigroups) of $\left(  m+n\right)  \times\left(
p+q\right)  $-supermatrices.
\end{problem}

First example of such new substructures are $\Gamma$-matrices considered below.

\begin{conjecture}
These new substructures lead to \textbf{new corresponding superoperators}
which are represented by one-parameter substructures of supermatrices.
\end{conjecture}

Therefore we first consider possible (not triangle) subsemigroups of supermatrices.

\section{Odd-reduced supermatrix semigroups}

In general, the odd-reduced matrices $M_{odd}$ do not form a semigroup, since
their multiplication is not closed in general $\mathbf{M}_{odd}\cdot
\mathbf{M}_{odd}\subset\mathbf{M}$. Nevertheless, some subset of
$\mathbf{M}_{odd}$ \textit{can form} a semigroup \cite{dup12}. That can happen
due to the existence of zero divisors in $\Lambda$, and so we have
$\mathbf{M}_{odd}\cdot\mathbf{M}_{odd}\cap\mathbf{M}_{odd}=\mathbf{M}%
_{odd}^{smg}\neq\emptyset$.

To find $\mathbf{M}_{odd}^{smg}$ we consider a $\left(  1+1\right)
\times\left(  1+1\right)  $ example. Let $\alpha,\beta\in\Gamma_{set}$, where
$\Gamma_{set}\subset\Lambda_{1}$. We denote $\mbox{\rm
Ann}\,\alpha\overset{def}{=}\left\{  \gamma\in\Lambda_{1}\,|\,\gamma
\cdot\alpha=0\right\}  $ and $\mbox{\rm
Ann}\,\Gamma_{set}=\mathrel{\mathop{\cap}\limits_{\alpha\in\Gamma
}}\mbox{\rm Ann}\,\alpha$ (here the intersection is crucial). Then we define
\textit{left} and \textit{right}\textsl{ }$\Gamma$\textit{-matrices}
\begin{equation}
\mathbf{M}_{odd(L)}^{\Gamma}\overset{def}{=}\left(
\begin{array}
[c]{cc}%
0 & \Gamma_{set}\\
\mbox{\rm
Ann}\,\Gamma_{set} & b
\end{array}
\right)  , \label{15c}%
\end{equation}%
\begin{equation}
\mathbf{M}_{odd(R)}^{\Gamma}\overset{def}{=}\left(
\begin{array}
[c]{cc}%
0 & \mbox{\rm Ann}\,\Gamma_{set}\\
\Gamma_{set} & b
\end{array}
\right)  . \label{15d}%
\end{equation}

\begin{proposition}
The $\Gamma$-matrices $\mathbf{M}_{odd(L,R)}^{\Gamma}\subset\mathbf{M}_{odd}$
form subsemigroups of $\frak{M}\left(  1|1\right)  $ under the standard
supermatrix multiplication, if $b\Gamma\subseteq\Gamma$.
\end{proposition}

\begin{definition}
$\Gamma$-semigroups\textsl{ }$\frak{M}_{odd(L,R)}^{\Gamma}\left(  1|1\right)
$ are subsemigroups of $\frak{M}\left(  1|1\right)  $ formed by the $\Gamma
$-matrices $\mathbf{M}_{odd(L,R)}^{\Gamma}$ under supermatrix multiplication.
\end{definition}

\begin{corollary}
The $\Gamma$-matrices are additional to triangle supermatrices substructures
which form semigroups.
\end{corollary}

Let us consider general square antitriangle $\left(  p+q\right)  \times\left(
p+q\right)  $-supermatrices (having even parity in notations of \cite{berezin}%
) of the form
\begin{equation}
M_{odd}^{p|q}\overset{def}{=}\left(
\begin{array}
[c]{cc}%
0_{p\times p} & \Gamma_{p\times q}\\
\Delta_{q\times p} & B_{q\times q}%
\end{array}
\right)  , \label{modd}%
\end{equation}
where ordinary matrix $B_{q\times q}$ consists of even elements and matrices
$\Gamma_{p\times q}$ and $\Delta_{q\times p}$ consist of odd elements
\cite{berezin,lei1} (we drop their indices below). The berezinian of
$M_{odd}^{p|q}$ can be obtained from the general formula by reduction and in
case of invertible $B$ (which is implied here) is (cf. (\ref{b0}))%
\begin{equation}
\operatorname*{Ber}M_{odd}^{p|q}=-\dfrac{\det\left(  \Gamma B^{-1}%
\Delta\right)  }{\det B}. \label{berm}%
\end{equation}

\begin{assertion}
A set of supermatrices $\mathbf{M}_{odd}^{p|q}$ form a semigroup
$\frak{M}_{odd}^{\Gamma}\left(  p|q\right)  $ of $\Gamma^{p|q}$-matrices, if
$\Gamma_{set}\Delta_{set}=0$, i.e. antidiagonal matrices are orthogonal, and
$\Gamma_{set}\mathbf{B}\subset\Gamma_{set}$, $\mathbf{B}\Delta_{set}%
\subset\Delta_{set}$.
\end{assertion}

\begin{proof}
Consider the product
\begin{equation}
M_{odd_{1}}^{p|q}M_{odd_{2}}^{p|q}=\left(
\begin{array}
[c]{cc}%
\Gamma_{1}\Delta_{2} & \Gamma_{1}B_{2}\\
B_{1}\Delta_{2} & B_{1}B_{2}+\Delta_{1}\Gamma_{2}%
\end{array}
\right)  \label{mm}%
\end{equation}
and observe the condition of vanishing even-even block, which gives
$\Gamma_{1}\Delta_{2}=0$, and other conditions follows obviously.
\end{proof}

From (\ref{mm}) it follows

\begin{corollary}
Two $\Gamma^{p|q}$-matrices satisfy the band relation $M_{1}M_{2}=M_{1}$, if
$\Gamma_{1}B_{2}=\Gamma_{1},$ $B_{1}\Delta_{2}=\Delta_{2},$ $B_{1}B_{2}%
+\Delta_{1}\Gamma_{2}=B_{1}$.
\end{corollary}

\begin{definition}
We call a set of $\Gamma^{p|q}$-matrices satisfying additional condition
$\Delta_{set}\Gamma_{set}=0$, a set of strong $\Gamma^{p|q}$-matrices.
\end{definition}

Strong $\Gamma^{p|q}$-matrices have some extra nice features and all
supermatrices considered below are of this class.

\begin{corollary}
Idempotent strong $\Gamma^{p|q}$-matrices are defined by relations $\Gamma
B=\Gamma,$ $B\Delta=\Delta,$ $B^{2}=B$.
\end{corollary}

The product of $n$ strong $\Gamma^{p|q}$-matrices $M_{i}$ has the following
form%
\begin{equation}
M_{1}M_{2}\ldots M_{n}=\left(
\begin{array}
[c]{cc}%
0 & \Gamma_{1}A_{n-1}B_{n}\\
B_{1}A_{n-1}\Delta_{n} & B_{1}A_{n-1}B_{n}%
\end{array}
\right)  , \label{mmn}%
\end{equation}
where $A_{n-1}=B_{2}B_{3}\ldots B_{n-1}$, and its berezinian is%
\begin{equation}
\operatorname*{Ber}\left(  M_{1}M_{2}\ldots M_{n}\right)  =-\dfrac{\det\left(
\Gamma_{1}A_{n-1}\Delta_{n}\right)  }{\det\left(  B_{1}A_{n-1}B_{n}\right)  }.
\label{bmn}%
\end{equation}

\section{One-even-parameter supermatrix idempotent semigroups}

Here we investigate one-even-parameter subsemigroups of $\Gamma$%
-semigroups\textsl{\ }and as a particular example for clearness of statements
consider $\frak{M}_{odd}\left(  1|1\right)  $, where all characteristic
features taking place in general $\left(  p+q\right)  \times\left(
p+q\right)  $ as well can be seen. These formulas will be applied for
establishing corresponding superoperator semigroup properties.

A simplest semigroup can be constructed from antidiagonal nilpotent
supermatrices of the shape
\begin{equation}
Y_{\alpha}\left(  t\right)  \overset{def}{=}\left(
\begin{array}
[c]{cc}%
0 & \alpha t\\
\alpha & 0
\end{array}
\right)  . \label{16y}%
\end{equation}
where $t\in\Lambda^{1|0}$ is an even parameter of the Grassmann algebra
$\Lambda$ which continuously ''numbers'' elements $Y_{\alpha}\left(  t\right)
$ and $\alpha\in\Lambda^{0|1}$ is a fixed odd element of $\Lambda$ which
''numbers'' the sets $\mathbf{Y}_{\alpha}=\mathrel{\mathop{\cup
}\limits_{t}}Y_{\alpha}\left(  t\right)  $.

\begin{definition}
The supermatrices $Y_{\alpha}\left(  t\right)  $ together with a null
supermatrix $Z\overset{def}{=}\left(
\begin{array}
[c]{cc}%
0 & 0\\
0 & 0
\end{array}
\right)  $ form a \textit{continuous null semigroup} $\frak{Z}_{\alpha}\left(
1|1\right)  =\left\{  \mathbf{Y}_{\alpha}\cup Z;\cdot\right\}  $ having the
null multiplication
\begin{equation}
Y_{\alpha}\left(  t\right)  Y_{\alpha}\left(  u\right)  =Z. \label{16yyz}%
\end{equation}
\end{definition}

\begin{assertion}
For any fixed $t\in\Lambda^{1|0}$ the set $\left\{  Y_{\alpha}\left(
t\right)  ,Z\right\}  $ is a 0-minimal ideal in $\frak{Z}_{\alpha}\left(
1|1\right)  $.
\end{assertion}

\begin{remark}
If we consider, for instance, a one-even-parameter odd-reduced supermatrix
$R_{\alpha}\left(  t\right)  =\left(
\begin{array}
[c]{cc}%
0 & \alpha\\
\alpha &  t
\end{array}
\right)  $, then multiplication of $R_{\alpha}\left(  t\right)  $ is not
closed since $R_{\alpha}\left(  t\right)  R_{\alpha}\left(  u\right)
=\left(
\begin{array}
[c]{cc}%
0 & \alpha u\\
\alpha t & tu
\end{array}
\right)  \notin\mathbf{R}_{\alpha}=\bigcup\limits_{t}R_{\alpha}\left(
t\right)  $. Any other possibility except ones considered below also do not
give closure of multiplication.
\end{remark}

Thus the only nontrivial closed systems of one-even-parameter odd-reduced
(antitriangle) $\left(  1+1\right)  \times\left(  1+1\right)  $ supermatrices
are $\mathbf{P}_{\alpha}=\mathrel{\mathop{\cup
}\limits_{t}}P_{\alpha}\left(  t\right)  $ where%
\begin{equation}
P_{\alpha}\left(  t\right)  \overset{def}{=}\left(
\begin{array}
[c]{cc}%
0 & \alpha t\\
\alpha & 1
\end{array}
\right)  \label{16}%
\end{equation}
and $\mathbf{Q}_{\alpha}=\mathrel{\mathop{\cup
}\limits_{t}}Q_{\alpha}\left(  u\right)  $ where%
\begin{equation}
Q_{\alpha}\left(  u\right)  \overset{def}{=}\left(
\begin{array}
[c]{cc}%
0 & \alpha\\
\alpha u & 1
\end{array}
\right)  . \label{16a}%
\end{equation}

First, we establish multiplication properties of supermatrices $P_{\alpha
}\left(  t\right)  $ and $Q_{\alpha}\left(  u\right)  $. Obviously, that they
are idempotent.

\begin{assertion}
Sets of idempotent supermatrices $\mathbf{P}_{\alpha}$ and $\mathbf{Q}%
_{\alpha}$ form left zero and right zero semigroups respectively with
multiplication
\begin{align}
P_{\alpha}\left(  t\right)  P_{\alpha}\left(  u\right)   &  =P_{\alpha}\left(
t\right)  ,\label{m111}\\
Q_{\alpha}\left(  t\right)  Q_{\alpha}\left(  u\right)   &  =Q_{\alpha}\left(
u\right)  . \label{m1q1}%
\end{align}
if and only if $\alpha^{2}=0$.
\end{assertion}

\begin{proof}
It simply follows from supermatrix multiplication law and general previous considerations.
\end{proof}

\begin{corollary}
The sets $\mathbf{P}_{\alpha}$ and $\mathbf{Q}_{\alpha}$ are rectangular bands
since%
\begin{align}
P_{\alpha}\left(  t\right)  P_{\alpha}\left(  u\right)  P_{\alpha}\left(
t\right)   &  =P_{\alpha}\left(  t\right)  ,\label{ppp1}\\
P_{\alpha}\left(  u\right)  P_{\alpha}\left(  t\right)  P_{\alpha}\left(
u\right)   &  =P_{\alpha}\left(  u\right)  \label{pp2}%
\end{align}
and%
\begin{align}
Q_{\alpha}\left(  u\right)  Q_{\alpha}\left(  t\right)  Q_{\alpha}\left(
u\right)   &  =Q_{\alpha}\left(  u\right)  ,\label{qqq1}\\
Q_{\alpha}\left(  t\right)  Q_{\alpha}\left(  u\right)  Q_{\alpha}\left(
t\right)   &  =Q_{\alpha}\left(  t\right)  \label{qqq2}%
\end{align}
with components $t=t_{0}+\operatorname*{Ann}\alpha$ and $u=u_{0}%
+\operatorname*{Ann}\alpha$ correspondingly.
\end{corollary}

They are orthogonal in sense of
\begin{equation}
Q_{\alpha}\left(  t\right)  P_{\alpha}\left(  u\right)  =E_{\alpha},
\label{qp}%
\end{equation}
where
\begin{equation}
E_{\alpha}\overset{def}{=}\left(
\begin{array}
[c]{cc}%
0 & \alpha\\
\alpha & 1
\end{array}
\right)  \label{e}%
\end{equation}
is a right ``unity'' and left ``zero'' in semigroup $\mathbf{P}_{\alpha}$,
because
\begin{equation}
P_{\alpha}\left(  t\right)  E_{\alpha}=P_{\alpha}\left(  t\right)
,\;\;E_{\alpha}P_{\alpha}\left(  t\right)  =E_{\alpha} \label{ep}%
\end{equation}
and a left ``unity'' and right ``zero'' in semigroup $\mathbf{Q}_{\alpha}$,
because
\begin{equation}
Q_{\alpha}\left(  t\right)  E_{\alpha}=E_{\alpha},\;\;E_{\alpha}Q_{\alpha
}\left(  t\right)  =Q_{\alpha}\left(  t\right)  . \label{eq}%
\end{equation}

It is important to note that
\begin{equation}
P_{\alpha}\left(  t=1\right)  =Q_{\alpha}\left(  t=1\right)  =E_{\alpha},
\label{pq1}%
\end{equation}
and so $\mathbf{P}_{\alpha}\cap\mathbf{Q}_{\alpha}=E_{\alpha}$. Therefore,
almost all properties of $\mathbf{P}_{\alpha}$ and $\mathbf{Q}_{\alpha}$ are
similar, and we will consider only one of them in what follows. For
generalized Green's relations and more detail properties of odd-reduced
supermatrices see \cite{dup15,dup-hab} .

\section{Odd-reduced supermatrix operator semigroups}

Let us consider a semigroup $\mathcal{P}$ of superoperators $\mathsf{P}\left(
t\right)  $ (see for general theory \cite{davies,goldstein,eng/nag})
represented by the one-even-parameter semigroup $\mathbf{P}_{\alpha}$ of
odd-reduced supermatrices $P_{\alpha}\left(  t\right)  $ (\ref{16}) which act
on $\left(  1|1\right)  $-dimensional superspace $\mathbb{R}^{1|1}$ as follows
$P_{\alpha}\left(  t\right)  \mathtt{X}$, where $\mathtt{X}=\left(
\begin{array}
[c]{c}%
x\\
\varkappa
\end{array}
\right)  \in\mathbb{R}^{1|1}$, where $x$ - even coordinate, $\varkappa$ - odd
coordinate $\left(  \varkappa^{2}=0\right)  $ having expansions (\ref{x1}) and
(\ref{x2}) respectively (see \textbf{Corollary \ref{cor}}). We have a
representation $\rho:\mathcal{P}\rightarrow\mathbf{P}_{\alpha}$ with
correspondence $\mathsf{P}\left(  t\right)  \rightarrow P_{\alpha}\left(
t\right)  $, but (as is usually made, e.g. \cite{eng/nag}) we identify space
of superoperators with the space of corresponding matrices (nevertheless, we
use here operator notations for convenience).

\begin{definition}
An odd-reduced ``dynamical'' system on $\mathbb{R}^{1|1}$ is defined by an
odd-reduced supermatrix-valued function $\mathsf{P}\left(  \cdot\right)
:\mathbb{R}_{+}\rightarrow\frak{M}_{odd}\left(  1|1\right)  $ and ``time
evolution'' of the state $\mathtt{X}\left(  0\right)  \in\mathbb{R}^{1|1}%
$given by the function $\mathtt{X}\left(  t\right)  :\mathbb{R}_{+}%
\rightarrow\mathbb{R}^{1|1}$, where
\begin{equation}
\mathtt{X}\left(  t\right)  =\mathsf{P}\left(  t\right)  \mathtt{X}\left(
0\right)  \label{xpx}%
\end{equation}
and can be called as orbit of $\mathtt{X}\left(  0\right)  $ under
$\mathsf{P}\left(  \cdot\right)  $.
\end{definition}

\begin{remark}
In general the definition, the continuity, the functional equation and most of
conclusions below hold valid also for $t\in\mathbb{R}^{1|0}$ (as e.g. in
\cite[p. 9]{eng/nag}) including ``nilpotent time'' directions (see
\textbf{Corollary \ref{cor}}).
\end{remark}

From (\ref{m111}) it follows that%
\begin{equation}
\mathsf{P}\left(  t\right)  \mathsf{P}\left(  s\right)  =\mathsf{P}\left(
t\right)  , \label{ps}%
\end{equation}
and so superoperators $\mathsf{P}\left(  t\right)  $ are idempotent. Also they
form a rectangular band, because of
\begin{align}
\mathsf{P}\left(  t\right)  \mathsf{P}\left(  s\right)  \mathsf{P}\left(
t\right)   &  =\mathsf{P}\left(  t\right)  ,\label{pppp1}\\
\mathsf{P}\left(  s\right)  \mathsf{P}\left(  t\right)  \mathsf{P}\left(
s\right)   &  =\mathsf{P}\left(  s\right)  . \label{pppp2}%
\end{align}

We observe that
\begin{equation}
\mathsf{P}\left(  0\right)  =\left(
\begin{array}
[c]{cc}%
0 & 0\\
\alpha & 1
\end{array}
\right)  \neq\mathsf{I}=\left(
\begin{array}
[c]{cc}%
1 & 0\\
0 & 1
\end{array}
\right)  , \label{pi}%
\end{equation}
as opposite to the standard case \cite{davies}. A ``generator'' $\mathsf{A}%
=\mathsf{P}^{\prime}\left(  t\right)  $ is
\begin{equation}
\mathsf{A}=\left(
\begin{array}
[c]{cc}%
0 & \alpha\\
0 & 0
\end{array}
\right)  , \label{a}%
\end{equation}
and so the standard definition of generator \cite{davies}
\begin{equation}
\mathsf{A}=\lim_{t\rightarrow0}\dfrac{\mathsf{P}\left(  t\right)
-\mathsf{P}\left(  0\right)  }{t}. \label{app}%
\end{equation}
holds and for difference we have the standard relation%
\begin{equation}
\mathsf{P}\left(  t\right)  -\mathsf{P}\left(  s\right)  =\mathsf{A}%
\cdot\left(  t-s\right)  . \label{ptu}%
\end{equation}
The following properties of the generator $\mathsf{A}$ take place
\begin{align}
\mathsf{P}\left(  t\right)  \mathsf{A}  &  =\mathsf{Z},\label{paz1}\\
\mathsf{AP}\left(  t\right)   &  =\mathsf{A}, \label{paz2}%
\end{align}
where ``zero operator'' $\mathsf{Z}$ is represented by the null supermatrix,
$\mathsf{A}^{2}=\mathsf{Z,}$ and therefore generator $\mathsf{A}$ is a
nilpotent of second degree.

From (\ref{app}) it follows that
\begin{equation}
\mathsf{P}\left(  t\right)  =\mathsf{P}\left(  0\right)  +\mathsf{A}\cdot t.
\label{pt}%
\end{equation}

\begin{definition}
We call operators which can be presented as a linear supermatrix function of
$t$ a $t$-linear superoperators.
\end{definition}

From (\ref{pt}) it follows that $\mathsf{P}\left(  t\right)  $ is a $t$-linear superoperator.

\begin{proposition}
Superoperators $\mathsf{P}\left(  t\right)  $ \textbf{cannot be presented} as
an exponent (as for the standard superoperator semigroups $\mathsf{T}\left(
t\right)  =e^{\mathsf{A}\cdot t}$ \cite{davies}).
\end{proposition}

\begin{proof}
In our case
\begin{equation}
\mathsf{T}\left(  t\right)  =e^{\mathsf{A}\cdot t}=\mathsf{I}+\mathsf{A}\cdot
t=\left(
\begin{array}
[c]{cc}%
1 & \alpha t\\
0 & 1
\end{array}
\right)  \notin\mathbf{P}_{\alpha}. \label{exp}%
\end{equation}
\end{proof}

\begin{remark}
Exponential superoperator $\mathsf{T}\left(  t\right)  =e^{\mathsf{A}\cdot t}$
is represented by even-reduced supermatrices $\mathsf{T}\left(  \cdot\right)
:\mathbb{R}_{+}\rightarrow\frak{M}_{even}\left(  1|1\right)  $ \cite{eng/nag},
but idempotent superoperator $\mathsf{P}\left(  t\right)  $ is represented by
odd-reduced supermatrices $\mathsf{P}\left(  \cdot\right)  :\mathbb{R}%
_{+}\rightarrow\frak{M}_{odd}\left(  1|1\right)  $ (see \textbf{Definition
\ref{even-odd}}).
\end{remark}

Nevertheless, the superoperator $\mathsf{P}\left(  t\right)  $ satisfies the
same linear differential equation
\begin{equation}
\mathsf{P}^{\prime}\left(  t\right)  =\mathsf{A}\cdot\mathsf{P}\left(
t\right)  \label{pap0}%
\end{equation}
as the standard exponential superoperator $\mathsf{T}\left(  t\right)  $ (the
initial value problem \cite{eng/nag})%
\begin{equation}
\mathsf{T}^{\prime}\left(  t\right)  =\mathsf{A}\cdot\mathsf{T}\left(
t\right)  . \label{tat}%
\end{equation}

That leads to the following

\begin{corollary}
In case initial state does not equal unity $\mathsf{P}\left(  0\right)
\neq\mathsf{I}$, there exists an additional class of solutions of the initial
value problem (\ref{pap0})-(\ref{tat}) among odd-reduced (antidiagonal)
idempotent $t$-linear (nonexponential) superoperators.
\end{corollary}

Let us compare behavior of superoperators $\mathsf{P}\left(  t\right)  $ and
$\mathsf{T}\left(  t\right)  $. First of all, their generators coincide
\begin{equation}
\mathsf{P}^{\prime}\left(  0\right)  =\mathsf{T}^{\prime}\left(  0\right)
=\mathsf{A}. \label{pta}%
\end{equation}
But powers of $\mathsf{P}\left(  t\right)  $ and $\mathsf{T}\left(  t\right)
$ are different $\mathsf{P}^{n}\left(  t\right)  =\mathsf{P}\left(  t\right)
$ and $\mathsf{T}^{n}\left(  t\right)  =\mathsf{T}\left(  nt\right)  $. In
their common actions the superoperator which is from the left transfers its
properties to the right hand side as follows
\begin{align}
\mathsf{T}^{n}\left(  t\right)  \mathsf{P}\left(  t\right)   &  =\mathsf{P}%
\left(  \left(  n+1\right)  t\right)  ,\label{tp1}\\
\mathsf{P}^{n}\left(  t\right)  \mathsf{T}\left(  t\right)   &  =\mathsf{P}%
\left(  t\right)  . \label{tp2}%
\end{align}

Their commutator is nonvanishing%
\begin{equation}
\left[  \mathsf{T}\left(  t\right)  \mathsf{P}\left(  s\right)  \right]
=\mathsf{P}^{\prime}\left(  0\right)  t=\mathsf{T}^{\prime}\left(  0\right)
t=\mathsf{A}t, \label{tpp}%
\end{equation}
which can be compared with the pure exponential commutator (for our case)
$\left[  \mathsf{T}\left(  t\right)  \mathsf{T}\left(  u\right)  \right]  =0$
and idempotent commutator%
\begin{equation}
\left[  \mathsf{P}\left(  t\right)  \mathsf{P}\left(  s\right)  \right]
=\mathsf{P}^{\prime}\left(  0\right)  \left(  t-s\right)  =\mathsf{A}\left(
t-s\right)  . \label{pppa}%
\end{equation}

\begin{assertion}
All superoperators $\mathsf{P}\left(  t\right)  $ and $\mathsf{T}\left(
t\right)  $ commute in case of ``nilpotent time'' and
\begin{equation}
t\in\operatorname*{Ann}\alpha. \label{ta}%
\end{equation}
\end{assertion}

\begin{remark}
The uniqueness theorem \cite[p. 3]{eng/nag} holds only for $\mathsf{T}\left(
t\right)  $, because the nonvanishing commutator $\left[  \mathsf{A}%
,\mathsf{P}\left(  t\right)  \right]  =\mathsf{A}\neq0$.
\end{remark}

\begin{corollary}
The superoperator $\mathsf{T}\left(  t\right)  $ is an inner inverse for
$\mathsf{P}\left(  t\right)  $, because of%
\begin{equation}
\mathsf{P}\left(  t\right)  \mathsf{T}\left(  t\right)  \mathsf{P}\left(
t\right)  =\mathsf{P}\left(  t\right)  , \label{ptp}%
\end{equation}
but it is not an outer inverse, because%
\begin{equation}
\mathsf{T}\left(  t\right)  \mathsf{P}\left(  t\right)  \mathsf{T}\left(
t\right)  =\mathsf{P}\left(  2t\right)  . \label{tpt}%
\end{equation}
\end{corollary}

Let us try to find a (possibly noninvertible) operator $\mathsf{U}$ which
connects exponential and idempotent superoperators $\mathsf{P}\left(
t\right)  $ and $\mathsf{T}\left(  t\right)  $.

\begin{assertion}
The ``semi-similarity'' relation%
\begin{equation}
\mathsf{T}\left(  t\right)  \mathsf{U}=\mathsf{UP}\left(  t\right)  \label{tu}%
\end{equation}
holds if
\begin{equation}
\mathsf{U}=\left(
\begin{array}
[c]{cc}%
\sigma\alpha & \sigma\\
0 & \rho\alpha
\end{array}
\right)  \label{u1}%
\end{equation}
which is noninvertible triangle and depends from two odd constants, and the
``adjoint'' relation%
\begin{equation}
\mathsf{U}^{\ast}\mathsf{T}\left(  t\right)  =\mathsf{P}\left(  t\right)
\mathsf{U}^{\ast} \label{ut}%
\end{equation}
holds if
\begin{equation}
\mathsf{U}^{\ast}=\left(
\begin{array}
[c]{cc}%
0 & \alpha vt\\
\alpha u & v
\end{array}
\right)  \label{u2}%
\end{equation}
which is also noninvertible antitriangle and depends from two even constants
and ``time''.
\end{assertion}

Note that $\mathsf{U}$ is nilpotent of third degree, since $\mathsf{U}%
^{2}=\sigma\rho\mathsf{A}$, but the ``adjoint'' superoperator is not nilpotent
at all if $v$ is not nilpotent.

Both $\mathsf{A}$ and $\mathsf{Z}$ behave as zeroes, but $\mathsf{Y}\left(
t\right)  $ (see (\ref{16y})) is a two-sided zero for $\mathsf{T}\left(
t\right)  $ only, since%
\begin{equation}
\mathsf{T}\left(  t\right)  \mathsf{Y}\left(  t\right)  =\mathsf{Y}\left(
t\right)  \mathsf{T}\left(  t\right)  =\mathsf{Y}\left(  t\right)  ,
\label{ty}%
\end{equation}
but
\begin{align}
\mathsf{P}\left(  t\right)  \mathsf{Y}\left(  t\right)   &  =\mathsf{Y}\left(
0\right)  ,\label{yp1}\\
\mathsf{Y}\left(  t\right)  \mathsf{P}\left(  t\right)   &  =\mathsf{A}t.
\label{yp2}%
\end{align}

If we add $\mathsf{A}$ and $\mathsf{Z}$ to superoperators $\mathsf{P}\left(
t\right)  $, then we obtain an extended odd-reduced noncommutative
superoperator semigroup $\mathcal{P}_{odd}=\bigcup\mathsf{P}\left(  t\right)
\bigcup\mathsf{A}\bigcup\mathsf{Z}$ with the following Cayley table (for
convenience we add $\mathsf{Y}\left(  t\right)  $ and $\mathsf{T}\left(
t\right)  $ as well)%

\begin{subequations}
\begin{equation}%
\begin{tabular}
[c]{|c||c|c|c|c||c|c|c|}\hline
$1\setminus2$ & $\mathsf{P}\left(  t\right)  $ & $\mathsf{P}\left(  s\right)
$ & $\mathsf{A}$ & $\mathsf{Z}$ & $\mathsf{Y}\left(  t\right)  $ &
$\mathsf{T}\left(  t\right)  $ & $\mathsf{T}\left(  s\right)  $\\\hline\hline
$\mathsf{P}\left(  t\right)  $ & $\mathsf{P}\left(  t\right)  $ &
$\mathsf{P}\left(  t\right)  $ & $\mathsf{Z}$ & $\mathsf{Z}$ & $\mathsf{P}%
\left(  t\right)  $ & $\mathsf{P}\left(  t\right)  $ & $\mathsf{P}\left(
t\right)  $\\\hline
$\mathsf{P}\left(  s\right)  $ & $\mathsf{P}\left(  s\right)  $ &
$\mathsf{P}\left(  s\right)  $ & $\mathsf{Z}$ & $\mathsf{Z}$ & $\mathsf{P}%
\left(  s\right)  $ & $\mathsf{P}\left(  s\right)  $ & $\mathsf{P}\left(
s\right)  $\\\hline
$\mathsf{A}$ & $\mathsf{A}$ & $\mathsf{A}$ & $\mathsf{Z}$ & $\mathsf{Z}$ &
$\mathsf{Z}$ & $\mathsf{A}$ & $\mathsf{A}$\\\hline
$\mathsf{Z}$ & $\mathsf{Z}$ & $\mathsf{Z}$ & $\mathsf{Z}$ & $\mathsf{Z}$ &
$\mathsf{Z}$ & $\mathsf{Z}$ & $\mathsf{Z}$\\\hline\hline
$\mathsf{Y}\left(  t\right)  $ & $\mathsf{A}t$ & $\mathsf{A}s$ & $\mathsf{Z}$%
& $\mathsf{Z}$ & $\mathsf{Z}$ & $\mathsf{Y}\left(  t\right)  $ &
$\mathsf{Y}\left(  t\right)  $\\\hline
$\mathsf{T}\left(  t\right)  $ & $\mathsf{P}\left(  2t\right)  $ &
$\mathsf{P}\left(  t+s\right)  $ & $\mathsf{A}$ & $\mathsf{Z}$ &
$\mathsf{Y}\left(  t\right)  $ & $\mathsf{T}\left(  2t\right)  $ &
$\mathsf{T}\left(  t+s\right)  $\\\hline
$\mathsf{T}\left(  s\right)  $ & $\mathsf{P}\left(  t+s\right)  $ &
$\mathsf{P}\left(  2s\right)  $ & $\mathsf{A}$ & $\mathsf{Z}$ & $\mathsf{Y}%
\left(  t\right)  $ & $\mathsf{T}\left(  t+s\right)  $ & $\mathsf{T}\left(
2s\right)  $\\\hline
\end{tabular}
\label{t}%
\end{equation}

It is easily seen that associativity in the left upper square holds, and so
the table (\ref{t}) is actually represents a semigroup of superoperators
$\mathcal{P}_{odd}$ (under supermatrix multiplication).

The analogs of the ``smoothing operator'' $\mathsf{V}\left(  t\right)  $
\cite{eng/nag} are%
\end{subequations}
\begin{align}
\mathsf{V}_{P}\left(  t\right)   &  =\int\limits_{0}^{t}\mathsf{P}\left(
s\right)  ds=\dfrac{t}{2}\left(  \mathsf{P}\left(  t\right)  +\mathsf{P}%
\left(  0\right)  \right)  =\left(
\begin{array}
[c]{cc}%
0 & \alpha\dfrac{t^{2}}{2}\\
\alpha t & t
\end{array}
\right)  ,\label{v1}\\
\mathsf{V}_{T}\left(  t\right)   &  =\int\limits_{0}^{t}\mathsf{T}\left(
s\right)  ds=\dfrac{t}{2}\left(  \mathsf{T}\left(  t\right)  +\mathsf{T}%
\left(  0\right)  \right)  =\left(
\begin{array}
[c]{cc}%
t & \alpha\dfrac{t^{2}}{2}\\
0 & t
\end{array}
\right)  . \label{v2}%
\end{align}

Let us consider the differential sequence of sets of superoperators
$\mathsf{P}\left(  t\right)  $%
\begin{equation}
\mathsf{S}_{n}\overset{\partial}{\rightarrow}\mathsf{S}_{n-1}\overset
{\partial}{\rightarrow}\ldots\mathsf{S}_{1}\overset{\partial}{\rightarrow
}\mathsf{S}_{0}\overset{\partial}{\rightarrow}\mathsf{A}\overset{\partial
}{\rightarrow}\mathsf{Z}, \label{ss2}%
\end{equation}
where $\partial=d/dt$ and
\begin{equation}
\mathsf{S}_{n}=\mathrel{\mathop{\bigcup}\limits_{t}}\dfrac{t^{n}}{n\left(
n-1\right)  \ldots1}\mathsf{P}\left(  \dfrac{t}{n+1}\right)  , \label{sp2}%
\end{equation}
and by definition
\begin{align}
\mathsf{S}_{0}  &  =\mathrel{\mathop{\bigcup}\limits_{t}}\mathsf{P}\left(
t\right)  ,\label{p2}\\
\mathsf{S}_{1}  &  =\mathrel{\mathop{\bigcup}\limits_{t}}\mathsf{V}_{P}\left(
t\right)  . \label{pv}%
\end{align}

Now we construct an analog of the standard operator semigroup functional
equation \cite{davies,eng/nag}
\begin{equation}
\mathsf{T}\left(  t+s\right)  =\mathsf{T}\left(  t\right)  \mathsf{T}\left(
s\right)  . \label{ttt}%
\end{equation}

Using the multiplication law (\ref{ps}) and manifest representation
(\ref{16}). for the idempotent superoperators $\mathsf{P}\left(  t\right)  $
we can formulate

\begin{definition}
The odd-reduced idempotent superoperators $\mathsf{P}\left(  t\right)  $
satisfy the following generalized functional equation
\begin{equation}
\mathsf{P}\left(  t+s\right)  =\mathsf{P}\left(  t\right)  \mathsf{P}\left(
s\right)  +\mathsf{N}\left(  t,s\right)  , \label{ppa}%
\end{equation}
where
\[
\mathsf{N}\left(  t,s\right)  =\mathsf{P}^{\prime}\left(  t\right)  s.
\]
\end{definition}

The presence of second term $\mathsf{N}\left(  t,s\right)  $ in the right hand
side of the generalized functional equation (\ref{ppa}) can be connected with
nonautonomous and deterministic properties of systems describing by it
\cite{eng/nag}. Indeed, from (\ref{xpx}) it follows that
\begin{align}
\mathtt{X}\left(  t+s\right)   &  =\mathsf{P}\left(  t+s\right)
\mathtt{X}\left(  0\right)  =\mathsf{P}\left(  t\right)  \mathsf{P}\left(
s\right)  \mathtt{X}\left(  0\right)  +\mathsf{P}^{\prime}\left(  t\right)
s\mathtt{X}\left(  0\right) \label{xxpp}\\
&  =\mathsf{P}\left(  t\right)  \mathtt{X}\left(  s\right)  +\mathsf{P}%
^{\prime}\left(  t\right)  s\mathtt{X}\left(  0\right)  \neq\mathsf{P}\left(
t\right)  \mathtt{X}\left(  s\right) \nonumber
\end{align}
as opposite to the always implied relation for exponential superoperators
$\mathsf{T}\left(  t\right)  $ (translational property \cite{davies,eng/nag})%
\begin{equation}
\mathsf{T}\left(  t\right)  \mathtt{X}\left(  s\right)  =\mathtt{X}\left(
t+s\right)  , \label{xts}%
\end{equation}
which follows from (\ref{ttt}). Instead of (\ref{xts}), using the band
property (\ref{ps}) we obtain%
\begin{equation}
\mathsf{P}\left(  t\right)  \mathtt{X}\left(  s\right)  =\mathtt{X}\left(
t\right)  , \label{xxp}%
\end{equation}
which can be called the ``moving time'' property.

\begin{problem}
Find a ``dynamical system'' with time evolution satisfying the ``moving time''
property (\ref{xxp}) instead of the translational property (\ref{xts}).
\end{problem}

\begin{assertion}
For ``nilpotent time'' satisfying (\ref{ta}) the generalized functional
equation (\ref{ppa}) coincides with the standard functional equation
(\ref{ttt}), and therefore the idempotent operators $\mathsf{P}\left(
t\right)  $ describe autonomous and deterministic ``dynamical'' system and
satisfy the translational property (\ref{xts}).
\end{assertion}

\begin{proof}
Follows from (\ref{ta}) and (\ref{xxpp}).
\end{proof}

\begin{problem}
\label{prob}Find all maps $\mathsf{P}\left(  \cdot\right)  :\mathbb{R}%
_{+}\rightarrow\frak{M}\left(  p|q\right)  $ satisfying the generalized
functional equation (\ref{ppa}).
\end{problem}

We turn to this problem later, and now consider some features of the Cauchy
problem for idempotent superoperators.

\section{Cauchy problem}

Let us consider an action (\ref{xpx}) of superoperator $\mathsf{P}\left(
t\right)  $ in superspace $\mathbb{R}^{1|1}$ as $\mathtt{X}\left(  t\right)
=\mathsf{P}\left(  t\right)  \mathtt{X}\left(  0\right)  $, where the initial
components are $\mathtt{X}\left(  0\right)  =\left(
\begin{array}
[c]{c}%
x_{0}\\
\varkappa_{0}%
\end{array}
\right)  $. From (\ref{xpx}) the evolution of the components has the form
\begin{equation}
\left(
\begin{array}
[c]{c}%
x\left(  t\right) \\
\varkappa\left(  t\right)
\end{array}
\right)  =\left(
\begin{array}
[c]{c}%
\alpha\varkappa_{0}t\\
\alpha x_{0}+\varkappa_{0}%
\end{array}
\right)  \label{xx}%
\end{equation}
which shows that superoperator $\mathsf{P}\left(  t\right)  $ does not lead to
time dependence of odd components. Then from (\ref{xx}) we see that
\begin{equation}
\mathtt{X}^{\prime}\left(  t\right)  =\left(
\begin{array}
[c]{c}%
\alpha\varkappa_{0}\\
0
\end{array}
\right)  =const. \label{x0}%
\end{equation}

This is in full agreement with an analog of the Cauchy problem for our case
\begin{equation}
\mathtt{X}^{\prime}\left(  t\right)  =\mathsf{A}\cdot\mathtt{X}\left(
t\right)  . \label{xax}%
\end{equation}

\begin{assertion}
The solution of the Cauchy problem (\ref{xax}) is given by (\ref{xpx}), but
the idempotent superoperator $\mathsf{P}\left(  t\right)  $ \textbf{can not be
presented} in exponential form as in the standard case \cite{davies}, but only
in the $t$-linear form $\mathsf{P}\left(  t\right)  =\mathsf{P}\left(
0\right)  +\mathsf{A}\cdot t\neq e^{\mathsf{A}\cdot t},$ as we have already
shown in (\ref{pt}).
\end{assertion}

This allows us to formulate

\begin{theorem}
In superspace the solution of the Cauchy initial problem with the same
generator $\mathsf{A}$ is two-fold and is given by two different type of superoperators:

\begin{enumerate}
\item  Exponential superoperator $\mathsf{T}\left(  t\right)  $ represented by
the even-reduced supermatrices;

\item  Idempotent $t$-linear superoperator $\mathsf{P}\left(  t\right)  $
represented by the odd-reduced supermatrices.
\end{enumerate}
\end{theorem}

For comparison the standard solution of the Cauchy problem (\ref{xax})
\[
\mathtt{X}\left(  t\right)  =\mathsf{T}\left(  t\right)  \mathtt{X}\left(
0\right)
\]
in components is%
\begin{equation}
\left(
\begin{array}
[c]{c}%
x\left(  t\right) \\
\varkappa\left(  t\right)
\end{array}
\right)  =\left(
\begin{array}
[c]{c}%
x_{0}+\alpha\varkappa_{0}t\\
\varkappa_{0}%
\end{array}
\right)  , \label{xxt}%
\end{equation}
which shows that the time evolution of even coordinate is also in nilpotent
even direction $\alpha\varkappa_{0}$ as in (\ref{xx}), but with addition of
initial (possibly nonilpotent) $x_{0}$, while odd coordinate is (another)
constant as well. That leads to

\begin{assertion}
``Even'' and ``odd'' evolutions coincide if even initial coordinate vanishes
$x_{0}=0$ or common starting point is pure odd $\mathtt{X}\left(  0\right)
=\left(
\begin{array}
[c]{c}%
0\\
\varkappa_{0}%
\end{array}
\right)  $.
\end{assertion}

A very much important formula is the condition of commutativity \cite{davies}
\begin{equation}
\left[  \mathsf{A},\mathsf{P}\left(  t\right)  \right]  \mathtt{X}\left(
t\right)  =\mathsf{A}\mathtt{X}\left(  t\right)  =\left(
\begin{array}
[c]{c}%
\alpha\varkappa\left(  t\right) \\
0
\end{array}
\right)  =0, \label{apx}%
\end{equation}
which satisfies, when $\alpha\cdot\varkappa\left(  t\right)  =0$, while in the
standard case the commutator $\left[  \mathsf{A},\mathsf{T}\left(  t\right)
\right]  \mathtt{X}\left(  t\right)  =0$, i.e. vanishes without any additional
conditions \cite{davies}.

\section{Superanalog of resolvent for exponential and idempotent superoperators}

For resolvents $\mathsf{R}_{P}\left(  z\right)  $ and $\mathsf{R}_{T}\left(
z\right)  $ we use analog the standard formula from \cite{davies} in the form
\begin{align}
\mathsf{R}_{P}\left(  z\right)   &  =\int\limits_{0}^{\infty}e^{-zt}%
\mathsf{P}\left(  t\right)  dt,\label{px}\\
\mathsf{R}_{T}\left(  z\right)   &  =\int\limits_{0}^{\infty}e^{-zt}%
\mathsf{T}\left(  t\right)  dt. \label{tx}%
\end{align}

Using the supermatrix representation (\ref{16}) we obtain
\begin{align}
\mathsf{R}_{P}\left(  z\right)   &  =\left(
\begin{array}
[c]{cc}%
0 & \dfrac{\alpha}{z^{2}}\\[12pt]%
\dfrac{\alpha}{z} & \dfrac{1}{z}%
\end{array}
\right)  ,\label{rz}\\
\mathsf{R}_{T}\left(  z\right)   &  =\left(
\begin{array}
[c]{cc}%
\dfrac{1}{z} & \dfrac{\alpha}{z^{2}}\\[12pt]%
0 & \dfrac{1}{z}%
\end{array}
\right)  . \label{rz1}%
\end{align}

We observe, that $\mathsf{R}_{T}\left(  z\right)  $ satisfies the standard
resolvent relation \cite{eng/nag}%
\begin{equation}
\mathsf{R}_{T}\left(  z\right)  -\mathsf{R}_{T}\left(  w\right)  =\left(
w-z\right)  \mathsf{R}_{T}\left(  z\right)  \mathsf{R}_{T}\left(  w\right)  ,
\label{rrt}%
\end{equation}
but its analog for $\mathsf{R}_{P}\left(  z\right)  $
\begin{equation}
\mathsf{R}_{P}\left(  z\right)  -\mathsf{R}_{P}\left(  w\right)  =\left(
w-z\right)  \mathsf{R}_{P}\left(  z\right)  \mathsf{R}_{P}\left(  w\right)
+\dfrac{w-z}{zw^{2}}\mathsf{A} \label{rra}%
\end{equation}
has additional term proportional to the generator $\mathsf{A}$.

\section{Properties of $t$-linear idempotent operators}

Here we consider properties of general $t$-linear (super)operators of the form%
\begin{equation}
\mathsf{K}\left(  t\right)  =\mathsf{K}_{0}+\mathsf{K}_{1}t, \label{k}%
\end{equation}
where $\mathsf{K}_{0}=\mathsf{K}\left(  0\right)  $ and $\mathsf{K}%
_{1}=\mathsf{K}^{\prime}\left(  0\right)  $ are constant (super)operators
represented by $\left(  n\times n\right)  $ matrices or $\left(  p+q\right)
\times$.$\left(  p+q\right)  $ supermatrices with $t$ (``time'') independent
entries. Obviously, that the generator of a general $t$-linear (super)operator
is
\begin{equation}
\mathsf{A}_{K}=\mathsf{K}^{\prime}\left(  0\right)  =\mathsf{K}_{1}.
\label{ak}%
\end{equation}

We will find system of equations for $\mathsf{K}_{0}$ and $\mathsf{K}_{1}$ for
some special cases appeared in above consideration.

\begin{assertion}
If a $t$-linear (super)operator $\mathsf{K}\left(  t\right)  $ satisfies the
band equation (\ref{ps})%
\begin{equation}
\mathsf{K}\left(  t\right)  \mathsf{K}\left(  s\right)  =\mathsf{K}\left(
t\right)  , \label{kk}%
\end{equation}
then it is idempotent and the constant component (super)operators
$\mathsf{K}_{0}$ and $\mathsf{K}_{1}$ satisfy the system of equations%
\begin{align}
\mathsf{K}_{0}^{2}  &  =\mathsf{K}_{0},\label{k1}\\
\mathsf{K}_{1}^{2}  &  =\mathsf{Z},\label{k2}\\
\mathsf{K}_{1}\mathsf{K}_{0}  &  =\mathsf{K}_{1},\label{k3}\\
\mathsf{K}_{0}\mathsf{K}_{1}  &  =\mathsf{Z}, \label{k4}%
\end{align}
from which it follows, that $\mathsf{K}_{0}$ is idempotent, $\mathsf{K}_{1}$
is nilpotent, and $\mathsf{K}_{1}$ is right divisor of zero and left zero for
$\mathsf{K}_{0}$.
\end{assertion}

For non-supersymmetric operators we have

\begin{corollary}
\label{cor1}The components of $t$-linear operator $\mathsf{K}\left(  t\right)
$ have the following properties: idempotent matrix $\mathsf{K}_{0}$ is similar
to an upper triangular matrix with $1$ on the main diagonal and nilpotent
matrix $\mathsf{K}_{1}$ is similar to an upper triangular matrix with $0$ on
the main diagonal \cite{okninski1,rad/ros}.
\end{corollary}

Comparing with the previous particular super case (\ref{pt}) we have
$\mathsf{K}_{0}=\mathsf{P}\left(  0\right)  $ and $\mathsf{K}_{1}%
=\mathsf{A}=\mathsf{P}^{\prime}\left(  0\right)  $.

\begin{remark}
In case of $\left(  p+q\right)  \times\left(  p+q\right)  $ supermatrices the
triangularization properties of \textbf{Corollary \ref{cor1}} do not hold
valid due to presence divisors of zero and nilpotents among entries (see
\textbf{Corollary \ref{cor}}), and so the inner structure of the component
supermatrices satisfying (\ref{k1})-(\ref{k4}) can be much different from the
standard non-supersymmetric case \cite{okninski1,rad/ros}.
\end{remark}

Let us consider the structure of $t$-linear operator $\mathsf{K}\left(
t\right)  $ satisfying the generalized functional equation (\ref{ppa}).

\begin{assertion}
If a $t$-linear (super)operator $\mathsf{K}\left(  t\right)  $ satisfies the
generalized functional equation%
\begin{equation}
\mathsf{K}\left(  t+s\right)  =\mathsf{K}\left(  t\right)  \mathsf{K}\left(
s\right)  +\mathsf{K}^{\prime}\left(  t\right)  s, \label{kka}%
\end{equation}
then its component (super)operators $\mathsf{K}_{0}$ and $\mathsf{K}_{1}$
satisfy the system of equations%
\begin{align}
\mathsf{K}_{0}^{2}  &  =\mathsf{K}_{0},\label{kf1}\\
\mathsf{K}_{1}^{2}  &  =\mathsf{Z},\label{kf2}\\
\mathsf{K}_{1}\mathsf{K}_{0}  &  =\mathsf{K}_{1},\label{kf3}\\
\mathsf{K}_{0}\mathsf{K}_{1}  &  =\mathsf{Z}, \label{kf4}%
\end{align}
\end{assertion}

We observe that the systems (\ref{k1})-(\ref{k4}) and (\ref{kf1})-(\ref{kf4})
are fully identical. It is important to observe the connection of the above
properties with the differential equation for $t$-linear (super)operator
$\mathsf{K}\left(  t\right)  $
\begin{equation}
\mathsf{K}^{\prime}\left(  t\right)  =\mathsf{A}_{K}\cdot\mathsf{K}\left(
t\right)  . \label{kak}%
\end{equation}

Using (\ref{ak}) we obtain the equation for components%
\begin{align}
\mathsf{K}_{1}^{2}  &  =\mathsf{Z},\label{kd1}\\
\mathsf{K}_{1}\mathsf{K}_{0}  &  =\mathsf{K}_{1}. \label{kd2}%
\end{align}

That leads to the following

\begin{theorem}
For any $t$-linear (super)operator $\mathsf{K}\left(  t\right)  =$
$\mathsf{K}_{0}+\mathsf{K}_{1}t$ the next statements are equivalent:

\begin{enumerate}
\item $\mathsf{K}\left(  t\right)  $ is idempotent and satisfies the band
equation (\ref{kk});

\item $\mathsf{K}\left(  t\right)  $ satisfies the generalized functional
equation (\ref{kka});

\item $\mathsf{K}\left(  t\right)  $ satisfies the differential equation
(\ref{kak}) and has idempotent time independent part $\mathsf{K}_{0}%
^{2}=\mathsf{K}_{0}$ which is orthogonal to its generator $\mathsf{K}%
_{0}\mathsf{A}=\mathsf{Z}$.
\end{enumerate}
\end{theorem}

\section{General $t$-power-type idempotent operators}

Let us consider idempotent (super)operators which depend from time by
power-type function, and so they have the form%
\begin{equation}
\mathsf{K}\left(  t\right)  =\sum\limits_{m=0}^{n}\mathsf{K}_{m}t^{m},
\label{km}%
\end{equation}
where $\mathsf{K}_{m}$ are $t$-independent (super)operators represented by
$\left(  n\times n\right)  $ matrices or $\left(  p+q\right)  \times$.$\left(
p+q\right)  $ supermatrices. This power-type dependence of is very much
important for super case, when supermatrix elements take value in Grassmann
algebra, and therefore can be nilpotent (see (\ref{x1})--(\ref{x2}) and
\textbf{Corollary \ref{cor}}).

We now start from the band property $\mathsf{K}\left(  t\right)
\mathsf{K}\left(  s\right)  =\mathsf{K}\left(  t\right)  $ and then find
analogs of the functional equation and differential equation for them.
Expanding the band property (\ref{kk}) in component we obtain $n$-dimensional
analog of (\ref{k1})-(\ref{k4}) as%
\begin{align}
\mathsf{K}_{0}^{2}  &  =\mathsf{K}_{0},\label{kn1}\\
\mathsf{K}_{i}^{2}  &  =\mathsf{Z},\;1\leq i\leq n,\label{kn2}\\
\mathsf{K}_{i}\mathsf{K}_{0}  &  =\mathsf{K}_{i},\;1\leq i\leq n,\label{kn3}\\
\mathsf{K}_{0}\mathsf{K}_{i}  &  =\mathsf{Z},\;1\leq i\leq n,\label{kn4}\\
\mathsf{K}_{i}\mathsf{K}_{j}  &  =\mathsf{Z},\;1\leq i,j\leq n,i\neq j.
\label{kn5}%
\end{align}

\begin{proposition}
The $n$-generalized functional equation for any $t$-power-type idempotent
(super)operators (\ref{km}) has the form%
\begin{equation}
\mathsf{K}\left(  t+s\right)  =\mathsf{K}\left(  t\right)  \mathsf{K}\left(
s\right)  +\mathsf{N}_{n}\left(  t,s\right)  , \label{kkn}%
\end{equation}
where%
\begin{equation}
\mathsf{N}_{n}\left(  t,s\right)  =\sum\limits_{m=1}^{n}\sum\limits_{l=m}%
^{n}\mathsf{K}_{l}\dfrac{l\left(  l-1\right)  \ldots\left(  l-m+1\right)
}{m!}s^{m}t^{l-m}. \label{n}%
\end{equation}
\end{proposition}

\begin{proof}
For the difference using the band property (\ref{kk}) we have $\mathsf{N}%
_{n}\left(  t,s\right)  =\mathsf{K}\left(  t+s\right)  -\mathsf{K}\left(
t\right)  \mathsf{K}\left(  s\right)  =\mathsf{K}\left(  t+s\right)
-\mathsf{K}\left(  t\right)  $. Then we expand in Taylor series around $t$ and
obtain $\mathsf{N}_{n}\left(  t,s\right)  =\sum\limits_{m=1}^{n}%
\mathsf{K}^{\left(  m\right)  }\left(  t\right)  \dfrac{s^{m}}{m!}$, where
$\mathsf{K}^{\left(  m\right)  }\left(  t\right)  $ denotes $n$-th derivative
which is a finite series for the power-type $\mathsf{K}\left(  t\right)  $
(\ref{km}).
\end{proof}

The differential equation for idempotent (super)operators coincide with the
standard initial value problem only for $t$-linear operators. In case of the
power-type operators (\ref{km}) we have

\begin{proposition}
The $n$-generalized differential equation for any $t$-power-type idempotent
(super)operators (\ref{km}) has the form%
\begin{equation}
\mathsf{K}^{\prime}\left(  t\right)  =\mathsf{A}_{K}\cdot\mathsf{K}\left(
t\right)  +\mathsf{U}_{n}\left(  t\right)  , \label{ka}%
\end{equation}
where%
\begin{equation}
\mathsf{U}_{n}\left(  t\right)  =\left\{
\begin{array}
[c]{ll}%
0 & n=1\\
\sum\limits_{m=2}^{n}m\mathsf{K}_{m}t^{m-1} & n\geq2
\end{array}
\right.  . \label{u}%
\end{equation}
\end{proposition}

\begin{proof}
To find the difference $\mathsf{U}_{n}\left(  t\right)  $ we use the expansion
(\ref{km}) and the band conditions for components (\ref{kn1})--(\ref{kn5}).
\end{proof}

\section{Conclusion}

In general one-parametric semigroups and corresponding superoperator
semigroups represented by antitriangle idempotent supermatrices and their
generalization for any dimensions $p,q,m,n$ have many unusual and nontrivial
properties \cite{dup12,dup15,dup-hab}. Here we considered only some of them
related to their connection with functional and differential equations. It
would be interesting to generalize the above constructions to higher
dimensions and to study continuity properties of the introduced idempotent
superoperators. These questions will be investigated elsewhere.

\bigskip

\begin{acknowledgement}
The author is grateful to Jan Okni\'{n}ski for valuable remarks and kind
hospitality at the Institute of Mathematics, Warsaw University, where this
work was begun. Also fruitful discussions with W. Dudek, A. Kelarev, G.
Kourinnoy, W. Marcinek and B. V. Novikov are greatly acknowledged.
\end{acknowledgement}

\end{document}